\documentclass{elsarticle}

\usepackage{amsfonts}
\usepackage{amsmath}
\usepackage{graphicx}
\usepackage[ansinew]{inputenc}

\newtheorem{teo}{Theorem}[section]
\newtheorem{prop}[teo]{Proposition}
\newtheorem{coro}[teo]{Corollary}
\newtheorem{lema}[teo]{Lemma}
\newtheorem{defi}[teo]{Definition}

\newtheorem{nota}{Remark}[section]

\def \e{\hskip 10pt }
\def \ecu{\hskip 5pt }
\def \Dem{{\bf Proof. }}

\def \sign{{\rm sign }}

\def \toro{\mathbb T}
\def \bbbz{\mathbb Z}
\def \bbbr{\mathbb R}

\def \bbbn{\mathbb N}
\newcommand{\re}{I \! \! R}

\def \num{${\mbox n}^{\underline{\mbox{\tiny o}}}$ }
\def \12{\frac{1}{2}}
\def \sinc{\mbox{\rm sinc}}

\def \dint {\int \! \! \! \int }

\def \hate {\hskip 2pt \hat{} \hskip 2pt}

\def \S{{\mathcal S}}

\def \M{{\mathcal M}}

\def \t2MB{\tilde{\tilde{{\mathcal MB}}}}

\def \ta{\hskip 2pt \widetilde{}\hskip 2pt }
\def \cqd{\vskip 20pt}
\def \supp{{\rm supp }}

\begin{document}

\begin{frontmatter}

\title{On boundedness of discrete multilinear singular integral operators}
\author[PV]{Paco Villarroya}
\ead{paco.villarroya@ed.ac.uk}
\fntext[PV]{The author has been partially supported by grants MTM 2005-08350-C03-03 and MTM2008-04594}
\address{School of Mathematics, University of Edinburgh, King's Buildings, EH10 3JZ, Edinburgh, UK}

\begin{abstract} Let $m(\xi,\eta)$ be a measurable locally bounded
function defined in $\mathbb R^2$. Let $1\leq p_1,q_1,p_2,q_2<\infty $ such that $p_i=1$ implies $q_i=\infty $. Let also $0<p_3,q_3<\infty $ and
$1/p=1/p_1+1/p_2-1/p_3$. We prove the following transference result: the operator
$$
{\mathcal C}_m(f,g)(x)=\int_{\bbbr} \int_{\bbbr} \hat f(\xi) \hat g(\eta)
m(\xi,\eta) e^{2\pi i x(\xi +\eta )}d\xi d\eta
$$
initially defined for integrable functions with compact Fourier support,
extends to a bounded bilinear operator from $L^{p_1,q_1}(\bbbr)\times
L^{p_2,q_2}(\bbbr)$ into $L^{p_3,q_3}(\bbbr)$ if and only if the family of operators
$$
{\mathcal D}_{\widetilde{m}_{t,p}} (a,b)(n)
=t^{\frac{1}{p}}\int_{-\12 }^{\12 }\int_{-\12 }^{\12 }P(\xi ) Q(\eta )
m(t\xi ,t\eta) e^{2\pi in(\xi +\eta  )}d\xi d\eta
$$
initially defined for finite sequences $a=(a_{k_{1}})_{k_{1}\in \bbbz }$, $b=(b_{k_{2}})_{k_{2}\in \bbbz }$,
where
$P(\xi )=\sum_{k_{1}\in \bbbz }a_{k_{1}}e^{-2\pi i k_{1}\xi }$ and
$Q(\eta )=\sum_{k_{2}\in \bbbz }b_{k_{2}}e^{-2\pi i k_{2}\eta }$, extend to bounded bilinear operators from
$l^{p_1,q_1}(\bbbz )\times l^{p_2,q_2}(\bbbz )$ into
$l^{p_3,q_3}(\bbbz )$ with norm bounded by uniform constant for all $t>0$.

We apply this result to prove boundedness of the discrete Bilinear Hilbert transforms and other related discrete multilinear
singular integrals including the endpoints.
\end{abstract}

\begin{keyword}
Multipliers \sep Transference \sep Lorentz spaces \sep Bilinear Hilbert transforms 
\MSC Primary 42A45, 42B20, 39A12; Secondary 45P05
\end{keyword}

\end{frontmatter}


\section{Introduction.}
Linear multiplier operators can be defined over a large variety of groups
in the following way:
given G a locally compact abelian group with Haar measure $\mu $
and  dual group $\hat{G}$, a measurable function
$m$ that takes values in $\hat{G}$ defines a
multiplier operator if for every $f\in L^{p}(G)$ there exists $g\in L^{p}(G)$ such that
${\mathcal F}g=m\cdot {\mathcal F}f$.
Here ${\mathcal F}$
is the continuous extension to $L^p(G)$ of the Fourier transform operator, initially defined in
$L^{1}(G)\cap L^{p}(G)$  as
${\mathcal F}f(\gamma )=\int_{G}f(x)\gamma (-x)d\mu (x)$ $\rm for$ $\rm every$ $\gamma \in \hat{G}$. The multiplier operator is then 
defined by $T_m(f)=g$.

This way, it is a very natural question to ask about the possibility of transferring
the boundedness properties of such operators when they are defined over two different groups. That is,
given a multiplier operator $T_m$
known to be bounded between spaces defined
over certain groups, let's say from $L^{p}(G_1)$ to $L^{q}(G_2)$,
we want to know when
the analogous operator $T_m'$ defined over different groups is also bounded between
similar type of spaces, let's say from $L^{p}(G_1')$ to $L^{q}(G_2')$.

The first transference methods for linear multipliers were given by K.
Deleeuw \cite{deleeuw} who showed that if $m$ is a bounded measurable
function which is pointwise limit of continuous functions then the linear operator
$$
T_m(f)(x)=\int_{\bbbr} \hat f(\xi)  m(\xi) e^{2\pi i x\xi }d\xi
$$
defined for $f\in S(\bbbr)$, extends boundedly to $L^p(\bbbr)$ with
$1\leq p<\infty $ if and only if
$$
\tilde{T}_{m_\varepsilon}(f)(\theta )
=\sum_{k\in \bbbz }\hat f(k)m(\varepsilon k)e^{2\pi i k\theta }
$$
defined for periodic functions $f$, extend to uniformly bounded operators on $L^p(\toro )$
for all $\varepsilon>0$.

Other type of linear transference theorems were given by P. Auscher and M.J. Carro (see \cite{aucarro})
on Lebesgue spaces between $\bbbr^n$ and
$\bbbz^n$. They proved that if $m$ is a measurable bounded function
then the operator
$$
T_m(f)(x)=\int_{\bbbr^{n}} \hat f(\xi)  m(\xi) e^{2\pi i x\xi }d\xi
$$
defined for $f\in S(\bbbr^{n})$, extend boundedly to $L^p(\bbbr^{n})$ with
$1\leq p<\infty $ if and only if
$$
\bar{T}_{m_\varepsilon}(a)(k)
=\int_{[-\frac{1}{2},\frac{1}{2}]^{n}}P(\xi )m(\varepsilon \xi )e^{2\pi ik\xi }d\xi
$$
for ${\displaystyle a=(a_{k})_{k\in \bbbz }}$ and
$P(\xi )=\sum_{k\in \bbbz }a_{k}e^{2\pi ik\xi }$
, extend to uniformly bounded operators on $l^p(\bbbz^{n})$ for all $\varepsilon>0$.

The interest for multilinear multipliers, which in the case of the real line can be defined as
$$C_m (f_1,f_2,\ldots ,f_n)(x)=\int_{\bbbr^n} \hat f_1(\xi_1)\ldots  \hat f_n(\xi_n)
m(\xi_1,\xi_2,\ldots ,\xi_n) e^{2\pi i x(\xi_1 +\xi_2+\ldots +\xi_n )}d\xi$$  for
$f_i\in{\mathcal S (\bbbr ) }$,
started in the seventies
with the works of R. Coifman and Y. Meyer. They proved (see \cite{CM1, CM2, CM3})
boundedness for multilinear multipliers
whose symbols $m$ have singularities at most at a single point.
At the end of the nineties, M. Lacey and C. Thiele
\cite{lath1, lath2} proved that the bilinear Hilbert transforms,
a family of bilinear multipliers for which
$m(\xi,\nu)=sign(\xi+\alpha\nu)$, are bounded multipliers from $L^{p_1}(\mathbb R)\times L^{p_2}(\mathbb R)$ into $L^{p_3}(\mathbb R)$
for $1<p_1,p_2\le\infty$, $p_3>2/3$ and $\alpha\in \bbbr\setminus\{ 0,1\}$.
Their paper was the first one with a proof of boundedness for multilinear multipliers whose symbols have singuralities spread over large sets.
This seminal work was quickly followed by many different extensions and generalizations. See the works by
Grafakos and Li \cite{grali1} and Li \cite{li}, by J.E.
Gilbert and A.R. Nahmod \cite {GN1, GN2}, by C. Muscalu, T. Tao and C.Thiele \cite{MTT, mutath2},
by M. Lacey \cite{lacey2} and by Grafakos, Tao and Terwilliger \cite{grataoter}.

Multilinear multiplier operators can also be defined over different groups 
and so the question of transference of boundedness properties also
applies to them. This way, D. Fan and S. Sato (see \cite{fansato})
proved the multilinear version of the transference between $\bbbr$ and
$\bbbz$, namely
that continuous functions $m(\xi,\eta)$ define multiplier operators of strong and weak
type $(p_1,p_2)$  on
$\mathbb R\times \mathbb R$ if and only if
$(m(\varepsilon k,\varepsilon k'))_{k,k'\in \mathbb Z}$ define a uniformly
bounded family of multipliers of strong and
weak type
$(p_1,p_2)$ on
$\bbbz\times\bbbz$.

Other references addressing the issue of transference of linear or multilinear multiplier operators through several different methods
are the following papers \cite{berblascagi}, \cite{blascargil},  \cite{caro}, \cite{carsor} and also the classic text \cite{cowe}.

The aim of the present paper is to get an extension of Auscher and Carro's result in the multilinear
setting for multipliers acting on Lorentz spaces which, in some sense, completes the Lorentz
transferences proven in \cite{blasvill} between $\bbbr $ and $\toro $. Unlike the linear case, in the multilinear setting
many interesting operators are known to be bounded in Lebesgue spaces with exponents below one and so these cases need to be
included in the transference results. This feature and the fact of
dealing with Lorentz norms are the main difficulties and novelties in the present work.
Although all results
hold true for multilinear multipliers in spaces of several variables $\mathbb R^n$ and $\mathbb Z^n$,
for the sake of simplicity in notation we restrict ourselves to bilinear operators with argument functions of one real variable.

We apply the transference results to prove $l^{p_1}(\mathbb Z)\times l^{p_2}(\mathbb Z)$ into $l^{p_3}(\mathbb Z)$
boundedness of the discrete Bilinear Hilbert transforms, defined for
any two finite sequences $a$, $b$ as
$$
{\mathcal H}_\alpha (a,b)(n)=\frac{1}{\pi }\sum_{k\in \mathbb Z,k\neq 0}a_{n-k}b_{n-\alpha k}\frac{1}{k}
$$
with $\alpha \in \mathbb Z\backslash \{ 0,1\} $.
This result has been previously proven by other methods (see \cite{blascargil}) when $\alpha =-1$
with the additional restriction $p_3\geq 1$. We manage to transfer boundedness for the whole family of operators 
even with exponents $p_3$ below one and also transfer weak boundedness at the endpoint
$p_3=2/3$.

We also notice that very minor changes can be done to tranfer
some operators with $x$-dependent symbols, that is operators whose representation is as follows
$$
{\mathcal C}_{m}(f,g)(x)
=\int_{\re }\hat{f}(\xi )\hat{g}(\eta )m(x,\xi ,\eta )e^{2\pi ix(\xi+\eta )}d\xi d\eta
$$
which allows to extend the applications to multilinear pseudodifferential operators.


\section{Lorentz spaces, interpolation and some notation.}
Let $X=(\Omega,\Sigma,\mu)$ be a $\sigma$-finite and complete measure space.
Given a measurable function $f:\Omega \rightarrow \mathbb C$ we shall denote its distribution function
by $\mu_f(\lambda) = \mu(\{ w \in \Omega : |f(w)| > \lambda\}$
and its noncreasing rearrangement function by
$f^*(t) = \inf\{\lambda>0 : \mu_f(\lambda) \leq t\}$.
The Lorentz space $L^{p,q}(X)$ is the family of all
measurable functions $f$ such that $\|f\|_{p,q} < \infty$, where
$$
\|f\|_{p,q} = \left\{
\begin{array}{ll}
\displaystyle{\left( \frac{q}{p} \int_0^\infty t^{\frac{q}{p}} f^*(t)^q \frac{dt}{t}\right)^\frac{1}{q},} &
0<p<\infty ,\ 0 < q < \infty , \\
\\
\displaystyle{\sup_{t> 0} t^{\frac{1}{p}} f^*(t) }  & 0 < p \leq \infty , \ q=\infty .
\end{array} \right.
$$
Lorentz spaces can be considered as a logarithmic refinement of Lebesgue spaces and actually when the exponents $p,q$ are equal, the related
Lorent space is the Lebesgue space $L^p(X)$.

We recall that simple functions are dense in
$L^{p,q}(X)$ for $q \neq \infty$ and that when $\Omega $ is finite dimensional and the measure 
is non-atomic the following dualities hold:
$(L^{p,q})^*(X) = L^{p',\infty}(X)$ for $1 \leq p < \infty $, $0<q\leq 1$;
$(L^{p,q})^*(X) = \{ 0\} $ for $p=1$, $1<q<\infty $ and
$(L^{p,q})^*(X) =L^{p',q'}(X)$ for $1 < p,q <\infty$.

On Lorentz spaces
the following version of Holder's inequality,
whose proof is due to O'Neil, holds:
if $0<p,p_1,p_2<\infty $ and $0<q,q_1,q_2\leq \infty $ obey $\frac{1}{p}=\frac{1}{p_1}+\frac{1}{p_2}$ and $\frac{1}{q}=\frac{1}{q_1}+\frac{1}{q_2}$ then
$$
\| fg\|_{L^{p,q}(X)}\leq C_{p_1,p_2,q_1,q_2}\| f\|_{L^{p_1,q_1}(X)}\| g\|_{L^{p_2,q_2}(X)}
$$
Moreover, when $p,q>1$, the following Minkowski's inequalities are satisfied
$$
\Big\| \sum_nf_n\Big\|_{L^{p,q}(X)}\leq C_{p,q}\sum_n\| f\|_{L^{p,q}(X)}
$$
$$
\Big\| \int fd\mu \Big\|_{L^{p,q}(X)}\leq C_{p,q}\int\| f\|_{L^{p,q}(X)}d\mu
$$
and also Young's inequality
$$
\| f*g\|_{L^{p,q}(X)}\leq C_{p_1,p_2,q}\| f\|_{L^{p_1,q}(X)}\| g\|_{L^{p_2,q}(X)}
$$
with $1<p,p_1,p_2,q$ and $\frac{1}{p}+1=\frac{1}{p_1}+\frac{1}{p_2}$.

It is well known the following representation of a Lorentz norm by the distribution function
$$
\|f\|_{p,\infty}= \sup_{\lambda>0} \lambda \mu_f(\lambda)^{1/p}
$$
We shall also use that for all $0<p,q<\infty$
\begin{eqnarray*}\label{equiv}
\|f\|_{p,q}&=&\left(
q \int_0^\infty \lambda^{q-1} \mu_f(\lambda)^{\frac{q}{p}}
d\lambda\right)^{\frac{1}{q}}
\end{eqnarray*}
which can be easily checked for simple functions. Two different changes of variables show that also
\begin{eqnarray*}\label{equiv}
\|f\|_{p,q}&=&\Big( \int_{0}^{\infty }
\mu_{f}(\lambda^{\frac{1}{q}})^{\frac{q}{p}}d\lambda \Big)^{\frac{1}{q}}
=\Big( \frac{q}{p}\int_{0}^{\infty }\lambda^{\frac{q}{p}-1}
\mu_{f}(\lambda^{\frac{1}{p}})^{\frac{q}{p}}d\lambda \Big)^{\frac{1}{q}}
\end{eqnarray*}
With these two expressions it is very easy to see that
\begin{equation*}\label{potencia}
\| f\|_{L^{p,q}(X)}\approx \| |f|^r\|_{L^{\frac{p}{r},\frac{q}{r}}(X)}^{\frac{1}{r}}
\end{equation*}
with actual equality in the case of the Lebesgue measure.

The reader is referred to \cite{BS}, \cite{berghlof},
\cite{gra}
or \cite{stewe}
for further information about Lorentz spaces.

\vskip 20pt

We end this section by setting up some tecnical notation and a definition.
For $0< p\leq \infty$, we denote the dilation, modulation and translation operators respectively by
$
D_{t}^{p}f(x)=t^{-\frac{1}{p}}f(t^{-1}x)$ (with the notation $D_t=D^\infty_t$), $M_{y}f(x)=f(x)e^{2\pi
iyx}
$ and
$
T_{y}f(x)=f(x-y)
$. They
satisfy the following symmetries through the Fourier transform:
$(T_{y}f)\hate{}=M_{-y}\hat{f}\ecu $,
$\ecu (M_{y}f)\hate{}=T_{y}\hat{f}\ecu $,
$\ecu (D_{t}^{p}f)\hate{}=\sign(t)D_{t^{-1}}^{p'}\hat{f} \ecu $ where, as usual, $p'$ stands for the conjugate
exponent of $p$. Moreover, we have dilation and translation invariance of Lorentz norms, that is
$$
\| D_{t}^{r}f\|_{L^{p,q}(X)}\approx |t|^{\frac{1}{p}-\frac{1}{r}}\| f\|_{L^{p,q}(X)}
\hskip 30pt
\| T_{y}f\|_{L^{p,q}(X)}=\| f\|_{L^{p,q}(X)}
$$

\vskip 10pt

We identify functions $f$
on $\toro $ and periodic functions on $\bbbr $  with period $1$ defined on $[-\12 ,\12 )$, that is
$f(x)=f(e^{2\pi ix})$  and
$
\int_{\toro }f(z)dm (z)=\int_{-\12 }^{\12 }f(t)dt
$. We give the following

\begin{defi}
We define the 1-periodization of a Schwartz function $f$ as
$$
\tilde{f}(x)=\sum_{k\in \bbbz }f(x+k)
$$
which, by the Poisson's summation formula, gives us the 1-periodic function whose Fourier coefficients
are given by the restriction of the Fourier transform of $f$ to $\bbbz $:
$$
\tilde{f}(x)=\sum_{k\in \bbbz }\hat{f}(k)e^{2\pi ix}
$$
\end{defi}

\section{Definition of the operators and statement of the main result.}
We define the operators whose boundedness properties we plan to transfer.
\begin{defi}
Let $m$ be a measurable locally bounded function defined in
$\bbbr^{2}$ and let ${\mathcal C}_{m}$ be the bilinear operator
$$
{\mathcal C}_{m}(f,g)(x)
=\dint_{\re^{2}}\hat{f}(\xi )\hat{g}(\eta )m(\xi ,\eta )
e^{2\pi i(\xi +\eta )x}d\xi d\eta
$$
initially defined for $f,g\in \S $ with compact Fourier support.

Let $0<p_{i},q_{i}\leq \infty $ with $i\in \{ 1,2,3\} $.
We say that $m$ is a
multiplier in $\bbbr $ if the operator can be extended to a bounded bilinear operator
from $L^{p_{1},q_{1}}(\bbbr )\times L^{p_{2},q_{2}}(\bbbr )$ to
$L^{p_{3},q_{3}}(\bbbr )$, that is if
$$
\| {\mathcal C}_{m}(f,g)\|_{L^{p_{3},q_{3}}(\bbbr )}\leq C\|f\|_{L^{p_{1},q_{1}}(\bbbr )}
\|g\|_{L^{p_{2},q_{2}}(\bbbr )}
$$
holds for all functions $f,g\in \S $ with compact Fourier support.
\end{defi}

We proved in \cite{vill} that a necessary condition for such a boundedness is that
if $p^{-1}=p_{1}^{-1}+p_{2}^{-1}-p_{3}^{-1}$ then $p\geq 1$, but this
constraint is no needed for a transference result.

We usually call multiplier both the function $m$ and the operator it defines,
but more often than not we reserve this notation for the operator and call
$m$ the symbol of the operator. If the symbol is an integrable
function the operator can be expressed via an integral kernel in the
following way
$$
{\mathcal C}_{m}(f,g)(x)=\dint_{\re^{2}}f(y)g(z)K(x-y,x-z)dydz
$$
where the kernel is $K=\check{m}$, that is the inverse Fourier transform of $m$.

Finally,
we notice that if $f,g\in \S $ then
${\mathcal C}_{m}(f,g)\in \S $ and
\begin{eqnarray}\label{tfoumulti}
{\mathcal C}_{m}(f,g)\hate{}(\nu )
=\int_{\re }\hat{f}(\xi )\hat{g}(\nu -\xi )m(\nu ,\nu -\xi )d\xi
\end{eqnarray}

\begin{defi}
Let $m$ be a measurable bounded function defined in
$\toro^{2}$.
We define the discrete bilinear operator ${\mathcal D}_m$, such that for all $n\in \bbbz $
$$
{\mathcal D}_{m}(a,b)(n)
=\int_{-\12 }^{\12 }\int_{-\12 }^{\12 }P(\xi )Q(\eta )
m(\xi ,\eta )e^{2\pi i(\xi +\eta )n}d\xi d\eta
$$
for finite sequences $a=(a_{k_{1}})_{k_{1}\in \bbbz }$,
$b=(b_{k_{2}})_{k_{2}\in \bbbz }$ where $P,Q$ are the trigonometric polynomials given by
$P(\xi )=\sum_{k_{1}\in \bbbz }a_{k_{1}}e^{-2\pi i k_{1}\xi }$ and
$Q(\eta )=\sum_{k_{2}\in \bbbz }b_{k_{2}}e^{-2\pi i k_{2}\eta }$.

We define the multipliers in
$\bbbz $ as those functions $m$ for which the previous operator
can be extended to a bounded bilinear operator from
$l^{p_{1},q_{1}}(\bbbz )\times l^{p_{2},q_{2}}(\bbbz )$ to
$l^{p_{3},q_{3}}(\bbbz )$.
\end{defi}
As before, if the symbol is a periodic integrable function then
the expresion of the
operator via an integral kernel is the following
$$
{\mathcal D}_{m}(a,b)(n)=\sum_{k_{1}\in \bbbz }\sum_{k_{2}\in \bbbz }
a_{k_{1}}b_{k_{2}}K(n-k_{1},n-k_{2})
$$
where now $K(n,l)=\widehat{D_{-1}m}(n,l)$ is the Fourier coefficient of $D_{-1}m$.

\vskip 10pt
Notice that in the transference between $\bbbr $ and $\toro $ (see \cite{blasvill}) the symbol
was discretized in order to define a periodic operator with
periodic argument functions.
Now, in the transference between $\bbbr $ and $\bbbz $, the symbol is periodized so that the operator and its argument
functions are discrete.

\vskip 20pt
Once all the necessary definitions have been established, we can summarize our main result in the following way.
Given a function $m$ defined in $\bbbr^{2}$, we prove that if $m$ defines a
bounded bilinear multiplier in $\bbbr $ so does in $\bbbz $ its periodization from $[-\12 ,\12 ]^2$,
$$
(m\chi_{[-\12 ,\12 ]^2})\ta (x,y)=\sum_{k_1,k_2\in \mathbb Z}(m\chi_{[-\12 ,\12 ]^2})((x,y)+(k_1,k_2))
$$
For the inverse transference, it is necessary that all periodizations of $m$ from growing intervals define a uniform
bounded family of bilinear multipliers in $\bbbz $ so that we recover that $m$ is also a bilinear multiplier in $\bbbr $.
The idea of the proof is to dilate the function $m$, to constraint this dilation to the interval
$[-\12 ,\12 ]^2$ and then to periodize. This way we get the family of symbols
$\widetilde{m}_{t,p}=(D_{t^{-1}}^{p}m\cdot \chi_{[-\12 ,\12 ]^2})\ta $ which need to be uniformly bounded for all $t>0$.

So the actual statement of the result is the following one:
\begin{teo}\label{teo-rtodeb} Let $m$ a locally bounded measurable function
defined in $\bbbr^{2}$. Let
$1\leq p_{i},q_{i}<\infty $ with $i=1,2$ such that $p_i=1$ implies $q_i=\infty $, $0< p_{3},q_{3}<\infty $ and
$p^{-1}=p_{1}^{-1}+p_{2}^{-1}-p_{3}^{-1}$.

Then $m$
is a multiplier in $\bbbr $ if and only if
$\{ \widetilde{m}_{t,p}\}_{t>0}$ defined in $\mathbb T^2$ is
a family of uniformly bounded multipliers in  $\bbbz $ for all
$t>0$. That is
$$
\| {\mathcal C}_{m}(f,g)\|_{L^{p_{3},q_{3}}(\bbbr )}\leq C\|f\|_{L^{p_{1},q_{1}}(\bbbr )}
\|g\|_{L^{p_{2},q_{2}}(\bbbr )}
$$
for all funtions $f,g\in \S $ if and only if
$$
\| {\mathcal D}_{\widetilde{m}_{t,p}}(a,b)\|_{L^{p_{3},q_3}(\bbbz )}\leq C
\|a\|_{L^{p_{1},q_1}(\bbbz )}\|b\|_{L^{p_{2},q_2}(\bbbz )}
$$
for all finite sequences $a=(a_{n})_{n}$, $b=(b_{n})_{n}$ and all $t>0$.
\end{teo}

\begin{nota}
Notice these two facts:

a) ${\mathcal C}_{m}$ is bounded if and only if ${\mathcal C}_{D_{t^{-1}}^{p}m}$ is bounded and both operators have the same norm.

b) if we define $\widetilde{m}_{t}=(D_{t^{-1}}m\cdot \chi_{[-\12 ,\12 ]^2})\ta $ then, 
the second inequality of the statement is equivalent to the fact that
$$
\| {\mathcal D}_{\widetilde{m}_t}(a,b)\|_{L^{p_{3},q_{3}}(\bbbz )}\leq Ct^{-\frac{1}{p}}
\|a\|_{L^{p_{1},q_{1}}(\bbbz )}\|b\|_{L^{p_{2},q_{2}}(\bbbz )}
$$
for all finite sequences $a=(a_{n})_{n}$, $b=(b_{n})_{n}$ and all $t>0$. But the equivalent boundedness of either ${\mathcal D}_{\widetilde{m}_{t,p}}$ or 
${\mathcal D}_{\widetilde{m}_{t}}$ for all $t>0$ are not equivalent to boundedness of
${\mathcal D}_{\widetilde{m}}$.
\end{nota}

\section{Equivalence between norms.}
In order to get the transference result, we will need two main ingredients: a relationship between the operators ${\mathcal C}_{m}$ and 
${\mathcal D}_{\tilde{m}}$ and
some equivalences between the norms of functions and sequences. This section is devouted to the latter.

We will first need
a relationship between the norm of a function and the norm of its discretization (restriction lemma) and then another relationship
between the norm of a sequence and the norm of certain function constructed with
such a sequence (extension lemma).
For that, we start by proving a relationship 'on average' and
then we move to a particular class of functions for which the equivalence of norms is somehow 'pointwise'.

\begin{lema}{\bf (Relationship on average)}\label{intesucedeb}
Let be $0<p,q\leq \infty $. For all $p_0<p<p_1$ there are constants $C_1,C_2>0$ such that 
for any $f\in \S $ and for each $n\in \mathbb Z$, if we define $a_n$ to be the function given by $a_{n}(u)=f(n+u)$ with $u\in [-\12 ,\12 )$, 
then
$$
C_{1}\Big\| \| a_{n}\|_{l^{p_0}(\mathbb T)}\Big\|_{l^{p,q}(\mathbb Z)}\leq
\| f\|_{L^{p,q}(\bbbr )}\leq C_2\Big\| \| a_{n}\|_{L^{p_1}(\mathbb T)}\Big\|_{l^{p,q}(\mathbb Z)}
$$
Here the implicit constants depend only on $p,q,p_0,p_1$.
\end{lema}

\Dem

Although this equivalence do not recover equality in the Lebesgue case, we show such a trivial identity
\begin{eqnarray*}
\| f\|_{L^{p}(\bbbr )}^{p}
=\sum_{n\in \bbbz }\int_{-\12 }^{\12 }|f(n+u)|^{p}du
&=& \Big\| \| a_{n}\|_{L^{p}(\mathbb T)}\Big\|_{l^{p}(\bbbz )}^p
\end{eqnarray*}
for all $p>0$ and similar for $p=\infty $. Now the result is deduced from this equality by interpolation.

We first assume that $\min(p,q)>1$.
For every $1\leq p_0<p<p_1$, we have by the Lebesgue case
$$
\| f\|_{L^{p_0}(\mathbb R)}
=\Big\| \| a_n\|_{L^{p_0}(\mathbb T)}\Big\|_{l^{p_0}(\mathbb Z )}
$$
while
$$
\| f\|_{L^{p_1}(\mathbb R)}=\Big\| \| a_n\|_{L^{p_1}(\mathbb T)}\Big\|_{l^{p_1}(\mathbb Z )}
\geq \Big\| \| a_n\|_{L^{p_0}(\mathbb T)}\Big\|_{l^{p_1}(\mathbb Z )}
$$

If we now define the sublinear operator
\begin{eqnarray*}
T_{p_0}:& L^p(\mathbb R)&\rightarrow l^p(\mathbb Z)\\
  & f &\hookrightarrow  (\| a_n\|_{L^{p_0}(\mathbb T)})_{n\in \mathbb Z}
\end{eqnarray*}
we have by the previous inequalities that
$
\| T_{p_0}(f)\|_{l^{p_i}(\mathbb Z)}\leq \| f\|_{L^{p_i}(\mathbb R)}
$
for $i=0,1$. So, by Marcinkiewicz's interpolation theorem 
$$
\Big\| \| a_n\|_{L^{p_0}(\mathbb T)}\Big\|_{l^{p,q}(\mathbb Z)}\leq C_{p_0,p,q}\| f\|_{L^{p,q}(\mathbb R)}
$$

To get the other inequality, we use duality, Holder's inequality and the previous case with $1\leq p_1'<p'$. For all $g\in {\mathcal S}(\mathbb R)$ with 
$\| g\|_{l^{p',q'}(\mathbb R)}=1$ and $b_n(u)=g(n+u)$,
$$
\Big| \int_{\mathbb R}f(x)g(x)dx\Big|=\Big|\sum_{n\in\bbbn }\int_{-\12 }^{\12 }f(n+u)g(n+u)du\Big|
\leq \sum_{n\in\bbbn }\| a_n\|_{L^{p_1}(\mathbb T )} \| b_n\|_{L^{p_1'}(\mathbb T )}
$$
$$
\leq \Big\|\| a_n\|_{L^{p_1}(\mathbb T )}\Big\|_{l^{p,q}(\mathbb Z )} \Big\| \| b_n\|_{L^{p_1'}(\mathbb T )}\Big\|_{l^{p',q'}(\bbbz )}
\leq \Big\| \| a_n\|_{L^{p_1}(\mathbb T )}\Big\|_{l^{p,q}(\mathbb Z )}C_{p_1',p',q'}\| g\|_{L^{p',q'}(\mathbb R)}
$$
that is
$$
\| f\|_{L^{p,q}(\mathbb R)}\leq C_{p_1',p',q'}\Big\| \| a_u\|_{L^{p_1}(\mathbb T )}\Big\|_{l^{p,q}(\mathbb Z )}
$$

When $\min(p,q)\leq 1$ we define $r=\min (p_0,q)-\epsilon <1$ and by previous cases
$$
\Big\| \| a_{n}\|_{L^{p_0}(\mathbb T)}\Big\|_{l^{p,q}(\mathbb Z)}
= \Big\| \| |a_{n}|^{r}\|_{L^{\frac{p_0}{r}}(\mathbb T )}\Big\|_{l^{\frac{p}{r},\frac{q}{r}}(\mathbb Z)}^{\frac{1}{r}}
\leq C\| |f|^{r}\|_{L^{\frac{p}{r},\frac{q}{r}}(\mathbb R)}^{\frac{1}{r}}
=C\| f\|_{L^{p,q}(\mathbb R)}
$$
while
$$
\| f\|_{L^{p,q}(\mathbb R)}=\| |f|^{r}\|_{L^{\frac{p}{r},\frac{q}{r}}(\mathbb R)}
\leq C \Big\| \| |a_n|^r\|_{L^{\frac{p_1}{r}}(\mathbb T)}\Big\|_{l^{\frac{p}{r},\frac{q}{r}}(\mathbb Z)}
=C \Big\| \| a_n \|_{L^{p_1}(\mathbb T)}\Big\|_{l^{p,q}(\mathbb Z)}
$$
This ends the proof.

\vskip 10pt


In general there is no relationship between the norm of a function and the norm as a sequence of its restriction to $\mathbb Z$. But
we will need an equivalence between both norms and this forces us to work with a class of functions for which
such an equivalence holds. This class will be 
the family functions of compact Fourier support.

\begin{defi} Given a distribution $u$ we say that $u\equiv 0$ in an open set $\Omega $ if
$(u,\varphi )=0$ for all $\varphi \in C^{\infty }_{0}(\Omega )$. Thus, we define the support
of a distribution $u$ as the complementary of the biggest open set $\Omega $ such that
$u\equiv 0$ in $\Omega$. 
\end{defi}
Notice that this set always exists. To prove it we first show that if a
distribution is null in an arbitrary family of open sets so is it in their union. Then the referred open set
$\Omega$ is simply the union of all open sets where $u$ is null.

\begin{defi}
A temperated distribution is called of exponential type (or of compact Fourier support) if its Fourier transform is
supported on a compact set. In particular, for $R>0$ we denote by $E_{R}$ the subspace of all temperated distributions whose Fourier transform is supported on the interval
$[-R,R]$.
\end{defi}

Two important properties of distributions in $E_{R}$ are the following theorems:
\begin{teo}{(Paley-Wiener's theorem).} Every distribution $f$
whose Fourier transform has compact suport in $[a,b]$
is the restriction to $\bbbr $ of an entire function
$F$ which besides satisfies the following bounds
$$
|F(x+it)|\leq C_{F}e^{-2at}
\e \e \e
|F(x-it)|\leq C_{F}e^{2bt}
$$
\end{teo}
This implies that $f$ is actually a function with full sense in each point and
moreover that $f\in C^{\infty }(\bbbr )$.

\begin{teo}{(Shannon's sampling theorem).} If $f\in L^{1}\cap E_{R}$ then
$$
f(x)=\sum_{n\in \bbbz }D_{2R}f(n)\hskip4pt \sinc{(2Rx-n)}
$$
where $\sinc (x)=\frac{1}{\pi x} \sin(\pi x)$.
\end{teo}
This implies that $f$ is totally determined by its restriction to $(2R)^{-1}\bbbz $.
Moreover, the result can be refined to get
$$
f(x)=\sum_{n\in \bbbz }D_{2R}f(n)\hskip4pt g(2Rx-n)
$$
where $g\in E_{1}\cap \S $ and $\hat{g}(x)=1$ for all $x\in [-\12 ,\12 ]$.

\vskip5pt
We also recall that ${\displaystyle \bigcup_{n\in \bbbn }E_{n}}$ is dense in $L^{p,q}$ for all possible
$0<p,q\leq \infty $.

\vskip 10pt

Roughly speaking, these results are two different ways of stating that functions of compact Fourier support
are essencially constant on intervals of lenght one and so
behave somehow like sequences: they are very smooth (without sudden spikes) and
fully determined by countable many samples.

The following lemma
is yet another way of expresing the same idea, this time by proving that
the norms of all discretizations of a function of compact Fourier support to lattices of lenght one
are controlled by the norm of the function itself.
In the case of $L^p$ norms with $p\geq 1$, this lemma is a classical result of entire function theory with many different known proofs.
The proof included here is a variation of that one in \cite{aucarro} that will be very useful
for our purposes.

\begin{lema}{\bf (Restriction)}\label{entereal}
Let $0<p,q\leq \infty $. Then, for every $p_0<p$ there exists
$C>0$ dependent of $p_0,p,q$ such that
for every $f\in E_{R}$ if we define the sequence $a=(a(n))_{n\in \bbbz }$ by $a(n)=f(n)$
then $a\in l^{p,q}(\bbbz )$ and
$$
\| a\|_{l^{p,q}(\bbbz )}\leq C\max (1,R^{\frac{1}{p_0}})\| f\|_{L^{p,q}(\bbbr )}
$$
\end{lema}
\Dem
Let $f\in E_R$ and $M=\max (1,R)$. Since $R\leq M$ implies $E_R\subset E_{M}$ we also have $f\in E_M$.

Let $\psi \in \S $ such that
$\hat{\psi }(x)=1$ if $x\in [-1,1]$ and with support in $[-2,2]$. Since $D_{M}\hat{\psi }\equiv 1$ in
$[-M,M]$ and has its support in $[-2M,2M]$, we have that
$\psi_M=D_{M^{-1}}^{1}\psi \in E_{2M}$. Moreover, for any $f\in E_{R}$ we obtain 
$\hat{f}=\hat{f}D_{M}\hat{\psi }$, that is, 
$f=f*\psi_M$.

Now for every $n\in \mathbb Z$ and $u\in [-\frac{1}{2},\frac{1}{2}]$ we define the function $a_{n}(u)=f(n+u)$
and then we have
$$
a(n)=f(n)=\int_{\re }f(n-x)\psi_M(x)dx
$$
$$
=\sum_{k\in \bbbz }\int_{-\12 +k}^{\12 +k}f(n-x)\psi_M(x)dx
=\sum_{k\in \bbbz }\int_{-\12 }^{\12 }a_{n}(-u-k)\psi_M(u+k)du
$$

We assume first that $\min (p,q)>1$ and choose $1\leq p_0<p$. By Holder's inequality
$$
|a(n)|\leq \sum_{k\in \bbbz }\Big| \int_{-\12 }^{\12 }T_{k}a_n(u) T_{-k}\psi_M(u)du\Big|
\leq \sum_{k\in \bbbz }\| T_{k}a_{n}\|_{L^{p_0}(\mathbb T)}\| T_{-k}\psi_M \|_{L^{p_0'}(\mathbb T)}
$$
By Minkowski's inequality, the average lemma \ref{intesucedeb} and the translation invariance of a Lorentz's norm we get 
$$
\| a\|_{l^{p,q}(\bbbz )}
\leq \sum_{k\in \bbbz }\Big\| \| T_{k}a_{n}\|_{L^{p_0}(\mathbb T)}\Big\|_{l^{p,q}(\bbbz )} \| T_{-k}\psi_M \|_{L^{p_0'}(\mathbb T)}
$$
$$
\leq C\sum_{k\in \bbbz }\| T_{k}f\|_{L^{p,q}(\mathbb R)}\| T_{-k}\psi_M \|_{L^{p_0'}(\mathbb T)}
$$
$$
=C \| f\|_{L^{p,q}(\mathbb R)}\sum_{k\in \bbbz }\| T_{-k}\psi_M \|_{L^{p_0'}(\mathbb T)}
$$
and we are finished once we prove 
\begin{eqnarray}\label{psi}
\sum_{k\in \bbbz }\| T_{-k}\psi_M \|_{L^{p_0'}(\mathbb T)}\leq CM^{\frac{1}{p_0}}
\end{eqnarray}
First,
$$
\| T_{-k}\psi_M \|_{L^{p_0'}(\mathbb T)}
=M\Big( \int_{-\12 }^{\12 }|\psi (M(x+k))|^{p'_0}dx \Big)^{\frac{1}{p'_0}}
=M^{\frac{1}{p_0}}\Big( \int_{-\frac{M}{2} }^{\frac{M}{2}}|\psi (x+Mk)|^{p'_0}dx \Big)^{\frac{1}{p'_0}}
$$
Now, for $k\neq 0$ we have
$
|\psi (x+Mk)|\leq C_2 (1+|x+Mk|^2)^{-1}\leq C(M^2k^2)^{-1}
$
since $|x|\leq M/2$ implies $|x+Mk|\geq M|k|-|x|\geq M|k|/2$. 
This way
\begin{eqnarray*}\label{controlaux}
\sum_{k\in \bbbz }\Big( \int_{-\frac{M}{2} }^{\frac{M}{2}}|\psi (x+Mk)|^{p'_0}dx \Big)^{\frac{1}{p'_0}}
\end{eqnarray*}
$$
\leq \Big( \int_{-\frac{M}{2} }^{\frac{M}{2}}|\psi (x)|^{p'_0}dx \Big)^{\frac{1}{p'_0}}
+\frac{C}{M^2}\sum_{k\neq 0 }\frac{1}{k^2}M^{\frac{1}{p'_0}}
\leq \| \psi \|_{L^{p'_0}(\mathbb R)}+C
$$

Notice that the thesis follows from the equality $f=f*\psi_R$, 
but we can get exactly the same result if $f$ satisfies $|f|\leq |f|*D_{R^{-1}}^{1}|\psi |$ even
though if $f\notin E_{R}$.

So, to prove the remaining cases we show first that though $f\in E_{R}$ does not imply $|f|^{r}\in E_R$ for $r\leq 1$, there always exists a constant $C>0$ independent of $f$ and $r$ such that
$$
|f|^{r}\leq C (|f|^{r}*D_{R^{-1}}^{1}|\psi |^{r})
$$
Let's see this point. Since $f\in E_R$ implies
$f=f*\psi_R $ we have
$$
|f(x)|\leq \int_{\re }|f(x-y)\psi_R (y)|dy
$$
The function inside the integral is again of exponential type since its Fourier transform is $(T_{-x}D_{-1}f\cdot \psi_R )\hate{}=
-M_{x}D_{-1}^{1}\hat{f}*D_{R}\hat{\psi }$ and then
$
\supp{(T_{-x}D_{-1}f\cdot \psi_R )\hate{}}
\subset \supp{(D_{-1}\hat{f})}+\supp{(D_{R}\hat{\psi })}\subset [-3R,3R]
$
. This implies
$
T_{-x}D_{-1}f\cdot \psi_R =(T_{-x}D_{-1}f\cdot \psi_R )
*\psi_{3R}
$,
that is
$$
|f(x-y)\psi_R (y)|
=\Big| \int_{\re }f(x-t)\psi_R (t)\psi_{3R}(y-t)dt\Big|
$$
$$
\leq \int_{\re }|f(x-t)\psi_R (t)|
3R\frac{C_{2}}{1+|3R(y-t)|^{2}}dt
$$
and so
$$
\frac{|f(x-y)\psi_R (y)|}
{CR\int_{\re }|f(x-t)\psi_R (t)|dt}\leq 1
$$
Now since $r\leq 1$
$$
\frac{|f(x-y)\psi_R(y)|^{r}}
{C^{r}R^{r}\Big( \int_{\re }|f(x-t)\psi_R (t)|dt\Big)^{r}}
\geq \frac{|f(x-y)\psi_R (y)|}
{CR\int_{\re }|f(x-t)\psi_R (t)|dt}
$$
and then by integration
$$
C^{1-r}R^{1-r}\int_{\re }|f(x-y)\psi_R (y)|^{r}dy
\geq \Big( \int_{\re }|f(x-t)\psi_R (t)|dt\Big)^{r}
\geq |f(x)|^{r}
$$
which finally gives us
$$
|f(x)|^{r}
\leq C^{1-r}\int_{\re }|f(x-y)|^{r}D_{R^{-1}}^{1}|\psi |^{r}(y)dy
\leq C(|f|^{r}*D_{R^{-1}}^{1}|\psi |^{r})(x)
$$
since $C$ can be chosen greater than one.

This is enough to conclude the statement. When $\min (p,q)\leq 1$, we take $\epsilon >0$ such that
$r=\min(p_0,q)-\epsilon <1$ and then although $|f|^{r}\notin E_{R}$
we still have the inequality
$
|f|^{r}\leq C (|f|^{r}*D_{R^{-1}}^{1}|\psi |^{r})
$.
So by the previous case with $1<\frac{p_0}{r}<\frac{p}{r}$
$$
\| a\|_{l^{p,q}(\bbbz )}
=\| |a|^{r}\|_{L^{\frac{p}{r},\frac{q}{r}}(\bbbz )}^{\frac{1}{r}}
$$
$$
\leq C\max (1,R)^{\frac{1}{\frac{p_0}{r}}\frac{1}{r}}\| |f|^{r}\|_{L^{\frac{p}{r },\frac{q}{r }}(\bbbr )}^{\frac{1}{r}}
=C\max (1,R^{\frac{1}{p_0}})\| f\|_{L^{p,q}(\bbbr )}
$$

\cqd

The following lemma states again that functions in $E_R$ behave as sequences, now in the sense that they satisfy a Young type inequality
when 'convoluted' with actual sequences.
\begin{lema}{\bf (Extension)}\label{realentelore}
Let $0<p,q\leq \infty $. Let $R>0$ and $\varphi \in E_{R}$.
Then for every $p_1>p$
there is a constant $C>0$ depending on $p_1,p,q$
such that for every sequence $a=(a_{n})$
$$
\Big \| \sum_{n\in \bbbz }a_{n}T_{n}\varphi \Big\|_{L^{p,q}(\bbbr )}
\leq C\| \varphi \|_{L^{s,q}(\mathbb R)}\max (1,R)^{\frac{1}{s}-\frac{1}{p_1}}\| a\|_{l^{p,q}(\bbbz )}
$$
for $s=\min(p,q,1)$.

Moreover, if $\hat{\varphi }$ is a linear multiplier in
$L^{p,q}(\bbbr )$ for $p,q>1$ we can substitute the $L^{s,q}$-norm of $\varphi $ by its norm as a multiplier
$\| \hat{\varphi }\|_{\M_{p,q}}$.
We call the function $\sum_{n\in \bbbz }a_{n}T_{n}\varphi $ the extension of the sequence $a$.
\end{lema}

\Dem
We take $\psi $ like in the previous lemma, $M=\max (1,R)$ and $\psi_M=D_{M^{-1}}^{1}\psi $.
We know that for all $\varphi \in E_{R}\subset E_{M}$ we have $\varphi =\varphi *\psi_{M}$
and so
$$
\sum_{k\in \bbbz }a_{k}T_{k}\varphi
=\Big( \sum_{k\in \bbbz }a_{k}T_{k}\psi_{M} \Big) *\varphi
$$
We assume first that $\min(p,q)>1$. Thus by Young's inequality on Lorentz spaces
$$
\Big\| \sum_{k\in \bbbz }a_{k}T_{k}\varphi \Big\|_{L^{p,q}(\bbbr )}
\leq \| \varphi \|_{L^{1,q}(\mathbb R)}\Big\| \sum_{k\in \bbbz }a_{k}T_{k}\psi_{M}\Big\|_{L^{p,q}(\bbbr )}
$$
or by the multiplier property
$$
\Big\| \sum_{k\in \bbbz }a_{k}T_{k}\varphi \Big\|_{L^{p,q}(\bbbr )}
\leq \| \hat{\varphi }\|_{\M_{p,q}(\mathbb R)}\Big\| \sum_{k\in \bbbz }a_{k}T_{k}\psi_{M}\Big\|_{L^{p,q}(\bbbr )}
$$

In any case, we have to deal with the same second factor. By the average lemma \ref{intesucedeb}
$$elsarticle, keywords
\Big\| \sum_{k\in \bbbz }a_{k}T_{k}\psi_{M}\Big\|_{L^{p,q}(\mathbb R)}
\leq C\Big\| \| (\sum_{k\in \bbbz }a_{k}T_{k}\psi_{M})_n\|_{L^{p_1}(\mathbb T)}\Big\|_{l^{p,q}(\mathbb Z )}
$$
and we transfer the dependence of the $n$-th term from the function $\psi_M$ to the sequence $a$ by a change of variables
$$
(\sum_{k\in \bbbz }a_{k}T_{k}\psi_{M})_n(u)=\sum_{k\in \bbbz }a_{k}\psi_{M}(u+n-k)
=\sum_{k\in \bbbz }a_{k+n}\psi_{M}(u-k)
$$
This way by Minkowski's inequality
$$
\| (\sum_{k\in \bbbz }a_{k}T_{k}\psi_{M})_n\|_{L^{p_1}(\mathbb T)}
=\Big( \int_{-\12 }^{\12 }|a_{k+n}T_k\psi_{M}(u)|^{p_1}\Big)^{\frac{1}{p_1}}
\leq \sum_{k\in \bbbz }|a_{k+n}|  \|T_{k}\psi_{M}\|_{L^{p_1}(\mathbb T)}
$$
Now by Minkowski's inequality again and the translation invariance of a Lorentz norm, we obtain
$$
\Big\| \sum_{k\in \bbbz }a_{k}T_{k}\psi_{M}\Big\|_{L^{p,q}(\mathbb R )}
\leq \sum_{k\in \bbbz }\| T_{-k}a\|_{l^{p,q}(\mathbb Z)} \| T_{k}\psi_{M}\|_{L^{p_1}(\mathbb T)}
$$
$$
\leq \| a\|_{l^{p,q}(\mathbb Z)} \sum_{k\in \bbbz }\|T_{k}\psi_{M}\|_{L^{p_1}(\mathbb T)}
\leq CM^{\frac{1}{p_1'}}\| a\|_{l^{p,q}(\mathbb Z)}
$$
where the last inequality follows from (\ref{controlaux}) in restriction lemma.

When $\min (p,q)\leq 1$ we take $r=\min(p,q)-\epsilon <1$ and proceed in a similar way we did in restriction lemma. Since we know that
$|\varphi |^{r}\leq C(|\varphi |^{r}*D_{M^{-1}}^{1}|\psi |^{{r}})$ we have 
$$
\Big| \sum_{k\in \bbbz }a_{k}T_{k}\varphi \Big|^{r}
\leq \sum_{k\in \bbbz }|a_{k}|^{r}T_{k}|\varphi |^{r}
$$
$$
\leq C\sum_{n\in \bbbz }|a_{k}|^{r}T_{k}(|\varphi |^{r}*D_{R^{-1}}^{1}|\psi |^{r})
=C\Big( \sum_{n\in \bbbz }|a_{k}|^{r}T_{k}D_{R^{-1}}^{1}|\psi |^{r}\Big)*|\varphi |^{r}
$$

Now,
since $\frac{p}{r},\frac{q}{r}>1$ we have by Young's inequality and the previous case
$$
\Big\| \sum_{k\in \bbbz }a_{k}T_{k}\varphi \Big\|_{L^{p,q}(\bbbr )}
=\Big\| \Big| \sum_{k\in \bbbz }a_{k}T_{k}\varphi \Big|^{r}\Big\|_{L^{\frac{p}{r},\frac{q}{r}}(\mathbb R )}^{\frac{1}{r}}
$$
$$
\leq C \Big\| \Big( \sum_{k\in \bbbz }|a_{k}|^{r}T_{k}D_{M^{-1}}^{1}
|\psi |^{r}\Big)*|\varphi |^{r}\Big\|_{L^{\frac{p}{r},\frac{q}{r}}(\mathbb R )}^{\frac{1}{r}}
$$
$$
\leq C\| |\varphi |^{r} \|_{L^{1,\frac{q}{r}}(\mathbb R)}^{\frac{1}{r}}
\Big\| \sum_{k\in \bbbz }|a_{k}|^{r}T_{k}D_{M^{-1}}^{1}|\psi |^{r}\Big\|_{L^{\frac{p}{r},\frac{q}{r}}(\bbbr )}^{\frac{1}{r}}
$$
$$
\leq C\| \varphi \|_{L^{r,q}(\mathbb R)}\| |a_n|^{r} \|_{l^{\frac{p}{r},\frac{q}{r}}(\bbbz )}^{\frac{1}{r}}
\Big( \sum_{k\in \bbbz }\| T_{k}D_{M^{-1}}^{1}|\psi |^{r}\|_{L^{\frac{p_1}{r}}(\toro )}\Big)^{\frac{1}{r}}
$$
$$
=C\| \varphi \|_{L^{r,q}(\mathbb R)}\| a_n\|_{l^{p,q}(\bbbz )}
\Big( \sum_{k\in \bbbz }\| T_{k}D_{M^{-1}}^{1}|\psi |^{r}\|_{L^{\frac{p_1}{r}}(\toro )}\Big)^{\frac{1}{r}}
$$
Then as we did before to prove (\ref{controlaux}) 
$$
\Big( \sum_{k\in \bbbz }\| T_{k}D_{M^{-1}}^{1}|\psi |^{r}\|_{L^{\frac{p_1}{r}}(\toro )}\Big)^{\frac{1}{r}}
=\Big( M^{1-\frac{r}{p_1}}\sum_{k\in \bbbz }\int_{-\frac{M}{2}}^{\frac{M}{2}} T_{kM}|\psi (u)|^{p_1}du\Big)^{\frac{1}{r}}
\leq CM^{\frac{1}{r}-\frac{1}{p_1}}
$$
and $r$ tends to $\min (p,q)$ when $\epsilon $ tends to zero.

\cqd

Now we can get the corollary we will need for the transference theorem:
\begin{coro}\label{reversolore}
Let $0<p,q\leq \infty $. Then there exists  $C>0$ dependent on $p,q$ such that
for every $f\in E_{R}$ with $R<1/2$ we have that if
$a=(a_{n})_{n\in \bbbz }$ is defined by $a_{n}=f(n)$ then
$$
\| f\|_{L^{p,q}(\bbbr )}\leq C\| a\|_{l^{p,q}(\bbbz )}
$$
\end{coro}
\Dem
For every
$f\in E_{R}\subset E_{1/2}$ we have
by Shannon's sampling theorem
$$
f(x)
=\sum_{n\in \bbbz }f(n)\varphi (x-n)
$$
with $\varphi \in E_{1}\cap \S$ such that $\hat{\varphi }(\xi )=1$ for $\xi \in [-\12 ,\12 ]$.
This way, by the extension lemma \ref{realentelore}
$$
\| f\|_{L^{p,q}(\bbbr )}
=\Big\| \sum_{n\in \bbbz }f(n)T_{n}\varphi \Big\|_{L^{p,q}(\bbbr )}
\leq C\|\varphi \|_{L^{s,q}(\bbbr )}\|a\|_{l^{p,q}(\bbbz )}
$$
with $s=\min(p,q,1)$.
\cqd

From restriction lemma \ref{entereal} and corollary \ref{reversolore} we finally obtain the equivalence between the norm of
a function in $E_{R}$ with $R<1/2$ and the norm of the sequence generated by its restriction to
$\mathbb Z $, that is

\begin{lema}{\bf (Equivalence)}\label{equivlore}
Let $0<p,q\leq \infty $.
If $f\in E_{R}$ with $R<1/2$ and for all $u\in [-\12 ,\12 ]$ we define the sequence $a_{u}(n)=f(n+u)$
then $\| f\|_{L^{p,q}(\bbbr )}\approx \|a_{u}\|_{L^{p,q}(\bbbz )}$ with bounds independent of $u$.
\end{lema}
\Dem
From $\widehat{T_{-u}f}=M_{u}\hat{f}$ we deduce that if $f\in E_{R}$ then so does
$T_{-u}f\in E_{R}$ for every
$u\in [-\12 ,\12 ]$. So,
\begin{enumerate}
\item[i)] by restriction lemma \ref{entereal} there exists $C>0$ such that
$$
\|a_{u}\|_{L^{p,q}(\bbbz )}\leq C\| T_{-u}f\|_{L^{p,q}(\bbbr )}
=C\| f\|_{L^{p,q}(\bbbr )}
$$
\item[ii)] by extension corollary \ref{reversolore} there exists $C'>0$ such that
$$
\| f\|_{L^{p,q}(\bbbr )}=\| T_{-u}f\|_{L^{p,q}(\bbbr )}\leq C'\| a_{u}\|_{l^{p,q}(\bbbz )}
$$
\end{enumerate}
\cqd

Notice that in $E_{R}$ with general $R$ we do not have the same equivalence.

%
%

\section{Proof of the transference theorem}
We prove now the transference theorem.

\Dem

($\Rightarrow $) Let $\hat{\varphi }\equiv \chi_{[-\12 ,\12 ]}$ which is known to be a linear multiplier 
for all $p,q>1$ and for $p=1$, $q=\infty $.

Given two finite sequences $a,b$ we fix $t>0$
and then
\begin{eqnarray*}\label{reloperators1}
{\mathcal D}_{\widetilde{m}_t}(a,b)(n)&=&\int_{-\12 }^{\12 }\int_{-\12 }^{\12 }
P(\xi )Q(\eta )m(t(\xi ,\eta ))e^{2\pi i(\xi +\eta )n}d\xi d\eta
\\
&=&\int_{\re }\int_{\re }
P(\xi )\widehat{\varphi }(\xi )\ecu
Q(\eta )\widehat{\varphi }(\eta )
\ecu m(t(\xi ,\eta ))\ecu e^{2\pi i(\xi +\eta )n}d\xi d\eta
\\
&=&{\mathcal C}_{D_{t^{-1}}m}(f,g)(n)
\end{eqnarray*}
where
$$
\hat{f}(\xi )=P(\xi )\widehat{\varphi }(\xi )=\sum_{k_{1}\in \bbbz }a_{k_{1}}e^{-2\pi ik_{1}\xi }\ecu
\widehat{\varphi }(\xi )
$$
that is
$$
f(x)
=\sum_{k_{1}\in \bbbz }a_{k_{1}}\int_{\re }\widehat{\varphi }(\xi )
e^{-2\pi i(k_{1}-x)\xi }d\xi
=\sum_{k_{1}\in \bbbz }a_{k_{1}}T_{k_{1}}\varphi (x)
$$
is the extension function of $a$. Since $\varphi \in E_{1}$, we have $f\in E_1$ and by the extension lemma \ref{realentelore}
$$
\| f\|_{L^{p_1,q_1}(\bbbr )}
=\big\| \sum_{k_{1}\in \bbbz }a_{k_{1}}T_{k_{1}}\varphi \big\|_{L^{p_1,q_1}(\bbbr )}
\leq \| \hat{\varphi }\|_{{\mathcal M}^{p_1,q_1}(\mathbb R)}\| a\|_{l^{p_1,q_1}(\bbbz )}=C\| a\|_{l^{p_1,q_1}(\bbbz )}
$$
for $p_1,q_1>1$ or $p_1=1$, $q_1=\infty $, and the same for $g$.

Moreover, by formula (\ref{tfoumulti}) we have
$$
{\mathcal C}_{{\mathcal D}_{t^{-1}}m}(f,g)\hate{}(\nu )
=\int_{\re }\hat{f}(\xi )\hat{g}(\nu -\xi )D_{t^{-1}}m(\nu ,\nu -\xi )d\xi
$$
and then $f,g\in E_{1}$ imply ${\mathcal C}_{D_{t^{-1}}m}(f,g)\in E_{2}$.
Thus we can also apply restriction lemma \ref{entereal} to get
$$
\| ({\mathcal C}_{D_{t^{-1}}m}(f,g)(n))_{n\in \mathbb Z}\|_{l^{p_{3},q_{3}}(\bbbz )}
\leq C\| {\mathcal C}_{D_{t^{-1}}m}(f,g)\|_{L^{p_{3},q_{3}}(\bbbr )}
$$

All this, the relationship between operators and the hypothesis give us the required transference:
$$
\| {\mathcal D}_{\widetilde{m}_t}(a,b)\|_{l^{p_{3},q_{3}}(\bbbz )}
=\| ({\mathcal C}_{D_{t^{-1}}m}(f,g)(n))_{n\in \mathbb Z}\|_{l^{p_{3},q_{3}}(\bbbz )}
$$
$$
\leq C\| {\mathcal C}_{D_{t^{-1}}m}(f,g)\|_{L^{p_{3},q_{3}}(\bbbr )}
\leq Ct^{-\frac{1}{p}}\|f\|_{L^{p_{1},q_{1}}(\bbbr )}\|g\|_{L^{p_{2},q_{2}}(\bbbr )}
$$
$$
\leq C t^{-\frac{1}{p}}\| a\|_{l^{p_{1},q_{1}}(\bbbz )}
\| b\|_{l^{p_{2},q_{2}}(\bbbz )}
$$

($\Leftarrow $) For the reverse implication we obtain first what essentially is the same relationship between the operators we saw before
but now starting with functions instead of sequences.

Let $f,g\in E_{R}$ and we take $k\in \bbbn $ with $k\geq 4R$ so that 
$\supp{\hat{f}},\supp{\hat{g}}\subset[-\frac{k}{2},\frac{k}{2}]$.
As usual, for each $x\in \bbbr $ we write $x=n+u$ with
$n\in \bbbz $ and $u\in [-\12 ,\12 ]$ and then we have
\begin{eqnarray*}
{\mathcal C}_{m}(f,g)(x)
&=&\int_{\re }\int_{\re }\hat{f}(\xi )\chi_{[-\frac{k}{2},\frac{k}{2}]}(\xi )\ecu
\hat{g}(\eta )\chi_{[-\frac{k}{2},\frac{k}{2}]}(\eta )\ecu m(\xi ,\eta )\ecu e^{2\pi i(\xi +\eta )x}d\xi d\eta
\\
&=&k^{2}\int_{\re }\int_{\re }\hat{f}(k\xi )\chi_{[-\frac{1}{2},\frac{1}{2}]}(\xi )\ecu
\hat{g}(k\eta )\chi_{[-\frac{1}{2},\frac{1}{2}]}(\eta )\ecu m(k(\xi ,\eta ))\ecu e^{2\pi ik(\xi +\eta )x}d\xi
d\eta \\
\end{eqnarray*}
or equivalently
\begin{eqnarray*}
{\mathcal C}_{m}(f,g)(k^{-1}(n+u))
&=&\int_{-\12 }^{\12 }\int_{-\12 }^{\12 }
k\hat{f}(k\xi )e^{2\pi iu\xi }\ecu k\hat{g}(k\eta )e^{2\pi iu\eta }\ecu
m(k(\xi ,\eta ))\ecu e^{2\pi i(\xi +\eta )n}\ecu d\xi d\eta \\
&=&{\mathcal D}_{\widetilde{m}_k}(a_{k,u},b_{k,u})(n)
\end{eqnarray*}
that is,
\begin{eqnarray}\label{reloperators2}
(D_{k}{\mathcal C}_{m}(f,g))_{u}(n)
={\mathcal D}_{\widetilde{m}_k}(a_{k,u},b_{k,u})(n)
\end{eqnarray}

We take $P_{k,u}$ to be the Fourier series whose $n$-th coefficient is $a_{k,u}(n)$. By definition $P_{k,u}$ is
the $1$-periodic function whose truncation to the interval $[-\frac{1}{2},\frac{1}{2}]$ has values $k\hat{f}(k\xi )e^{2\pi iu\xi }$.
This way 
$$
a_{k,u}(n)=\int_{-\12 }^{\12 }P_{k,u}(\xi )e^{2\pi in\xi }d\xi
$$
$$
=\int_{-\12 }^{\12 }k\hat{f}(k\xi )e^{2\pi iu\xi }e^{2\pi in\xi }d\xi
=\int_{-\frac{k}{2}}^{\frac{k}{2}}\hat{f}(\xi )e^{2\pi ik^{-1}(n+u)\xi }d\xi
$$
and since $\supp{\hat{f}}\subset [-\frac{k}{2},\frac{k}{2}]$, by the inversion formula we obtain
$$
a_{k,u}(n)=\int_{\re }\hat{f}(\xi )e^{2\pi ik^{-1}(n+u)\xi }d\xi
=f(k^{-1}(n+u))
=T_{-u}D_{k}f(n)
$$

Now $f,g\in E_{\frac{k}{2}}$ imply
$D_{k}f,D_{k}g\in E_{\frac{1}{2}}$ and so, by restriction lemma \ref{entereal},
$$
\| a_{k,u}\|_{l^{p_{1},q_{1}}(\bbbz )}
\leq C\| T_{-u}D_{k}f\|_{L^{p_{1},q_{1}}(\bbbr )}=Ck^{\frac{1}{p_{1}}}\| f\|_{L^{p_{1},q_{1}}(\bbbr )}
$$
and the same for $b_{k,u}$.

We know that $f,g\in E_R$ imply ${\mathcal C}_m(f,g)\in E_{2R}$ and so $D_k{\mathcal C}_m(f,g)\in E_{2R/k}$ with $2R/k<1/2$.
Then by corollary \ref{reversolore}, the relationship between operators given by (\ref{reloperators2}), the hypothesis and the homogeneity relation 
$p_1^{-1}+p_2^{-1}=p_3^{-1}+p^{-1}$ we finally have
$$
\| {\mathcal C}_{m}(f,g)\|_{L^{p_{3},q_{3}}(\bbbr )}
=k^{-\frac{1}{p_{3}}}\| D_{k}{\mathcal C}_{m}(f,g)\|_{L^{p_{3},q_{3}}(\bbbr )}
$$
$$
\leq Ck^{-\frac{1}{p_{3}}}\| (D_{k}{\mathcal C}_{m}(f,g))_{u}\|_{l^{p_{3},q_{3}}(\mathbb Z )}
=Ck^{-\frac{1}{p_{3}}}
\| {\mathcal D}_{\widetilde{m}_k}(a_{k,u},b_{k,u})\|_{l^{p_{3},q_{3}}(\bbbz )}
$$
$$
\leq Ck^{-\frac{1}{p_{3}}}k^{-\frac{1}{p}}\| a_{k,u}\|_{L^{p_{1},q_{1}}(\bbbz )}
\| b_{k,u}\|_{L^{p_{2},q_{2}}(\bbbz )}
$$
$$
\leq Ck^{-(\frac{1}{p_{3}}+\frac{1}{p})}k^{\frac{1}{p_{1}}}\| f\|_{L^{p_{1},q_{1}}(\bbbr )}
k^{\frac{1}{p_{2}}}\| g\|_{L^{p_{2},q_{2}}(\bbbr )}
$$
$$
=C\| f\|_{L^{p_{1},q_{1}}(\bbbr )}\| g\|_{L^{p_{2},q_{2}}(\bbbr )}
$$

\section{Applications}
As said in the introduction, Lacey and Thiele (see \cite{lath1} and \cite{lath2}) showed boundedness for singular integrals operators
with singularities spread over large sets when
they proved boundedness of Bilinear Hilbert transforms. Their result can be stated as follows:
\begin{teo}\label{laceythiele}
Let the Bilinear Hilbert transform be
$$
H_\alpha (f,g)=\frac{1}{\pi }\lim_{\epsilon \rightarrow 0}\int_{|t|>\epsilon }f(x-t)g(x-\alpha t)\frac{1}{t}dx
$$
initially defined for Schwartz functions $f$, $g$.

Then, for every $\alpha\in \bbbr \setminus\{ 0,1\}$ and $1<p_1,p_2\le\infty$, $2/3<p_3$ such that
$1/p_1+1/p_2=1/p_3$ there exists
$C=C(\alpha,p_1,p_2)>0$ for which
$$
\| H_{\alpha}(f,g)\|_{p_3}\le C\|f\|_{p_1}\|g\|_{p_2}
$$
for all $f,g$ in the Schwartz class.
\end{teo}

Short after them, Bilyk and Grafakos proved in \cite{bilykgrafakos} several distributional estimates at the endpoints
either when one of the exponents $p_1$ or $p_2$ are equal to one or
when $p_3=2/3$.
Their results imply in particular the following
\begin{teo}\label{bilykgrafa}
Then, for every $\alpha\in \bbbr \setminus\{ 0,1\}$ and $1\leq p_1,p_2\le 2$, $2/3=p_3$ such that
$1/p_1+1/p_2=1/p_3$ there exists
$C=C(\alpha,p_1,p_2)>0$ for which
$$
\| H_{\alpha}(f,g)\|_{2/3,\infty }\le C\|f\|_{p_1}\|g\|_{p_2}
$$
for all $f,g$ in the Schwartz class.
\end{teo}

Let's see now how the transference result \ref{teo-rtodeb} allows us to prove the analogous discrete versions of such results. As said, boundedness
of discrete Bilinear Hilbert transform has been previously proven in \cite{blascargil} for $\alpha =-1$ and $p_3\geq1$. Now we extend the
result to the whole family including exponents $p_3$ below one
and prove weak type estimates at the endpoint $p_3=2/3$.
\begin{prop}\label{discretelaceythiele}
Let $\alpha \in \mathbb Z\backslash \{ 0,1\} $ and
$$
{\mathcal H}_\alpha (a,b)(n)=\frac{1}{\pi }\sum_{k\in \mathbb Z\backslash \{ 0\} }a_{n-k}b_{n-\alpha k}\frac{1}{k}
$$
initially defined for finite sequences $a$, $b$.

Then, for $1<p_1,p_2\le\infty$ and $2/3<p_3$ such that
$1/p_1+1/p_2=1/p_3$ there exists
$C=C(\alpha,p_1,p_2)>0$ for which
$$
\|{\mathcal H}_{\alpha}(f,g)\|_{l^{p_3}(\mathbb Z)}\leq C\|a\|_{l^{p_1}(\mathbb Z)}\|b\|_{l^{p_2}(\mathbb Z)}
$$
for all $a,b$ finite sequences. Moreover, when $1<p_1,p_2\leq 2$ and $p_3=2/3$ there exists
$C=C(\alpha,p_1,p_2)>0$ such that
$$
\|{\mathcal H}_{\alpha}(f,g)\|_{l^{2/3,\infty }(\mathbb Z)}\leq C\|a\|_{l^{p_1}(\mathbb Z)}\|b\|_{l^{p_2}(\mathbb Z)}
$$
for all $a,b$ finite sequences.

\end{prop}

\Dem
The representation of the Bilinear Hilbert transforms via Fourier transform is the following
$$
H_\alpha (f,g)=-i\int_{\mathbb R}\int_{\mathbb R}\hat{f}(\xi)\hat{g}(\eta )\sign (\xi +\alpha \eta )e^{2\pi i(\xi +\eta )x}d\xi d\eta
$$
that is $H_\alpha ={\mathcal C}_{m_\alpha }$ with $m_{\alpha }(\xi ,\eta )=-i\sign (\xi +\alpha \eta )$.

So theorems \ref{laceythiele} and \ref{bilykgrafa} and the transference result \ref{teo-rtodeb} prove boundedness also for the corresponding
${\mathcal D}_{\tilde{m}_{\alpha ,t,p}}$, where
$\widetilde{m}_{\alpha, t,p}=(D_{t^{-1}}^{p}m_{\alpha }\cdot \chi_{[-\12 ,\12 ]^2})\ta $.
But in this case $p^{-1}=p_1^{-1}+p_2^{-1}-p_3^{-1}=0$
and moreover $D_{t^{-1}}m_{\alpha }$ is independent of t. This implies that $\widetilde{m}_{\alpha, t,p}=(m_{\alpha }\cdot \chi_{[-\12 ,\12 ]^2})\ta $
and so we have boundedness
for ${\mathcal D}_{\alpha }:={\mathcal D}_{(m_{\alpha }\cdot \chi_{[-\12 ,\12 ]^2})\ta }$ when $\alpha \in \mathbb R\backslash \{ 0,1\} $.

We now calculate the operator ${\mathcal D}_{\alpha }$ in order to see its relationship with the discrete Bilinear Hilbert transform.
\begin{eqnarray*}
{\mathcal D}_{\alpha }(a,b)(n)&=&-i\int_{-\12 }^{\12 }\int_{-\12 }^{\12 }
P(\xi )Q(\eta )m_{\alpha }(\xi ,\eta )e^{2\pi i(\xi +\eta )n}d\xi d\eta \\
&=&-i\sum_{k_{1},k_{2}\in \bbbz }a_{k_{1}}b_{k_{2}}\int_{-\12 }^{\12 }\int_{-\12 }^{\12 }\sign (\xi +\alpha \eta )e^{2\pi i((n-k_{1})\xi +(n-k_{2})\eta })d\xi d\eta \\
&=&-i\sum_{k_{1},k_{2}\in \bbbz }a_{k_{1}}b_{k_{2}}c_{\alpha }(n-k_1,n-k_2)
\end{eqnarray*}
and we compute
$
c_\alpha (r,s)=\int_{-\12 }^{\12 }\int_{-\12 }^{\12 }\sign (\xi +\alpha \eta )
e^{2\pi i(r\xi +s\eta )}d\xi d\eta
$
for $|\alpha |\leq 1$. An easy calculation gives for $r\neq 0$,
\begin{eqnarray*}
c_\alpha (r,s)&=&-\frac{1}{\pi ir}\Big( \sinc(-\alpha r+s)-\cos(\pi r)\, \sinc{(s)}\Big) \\
&=&-\frac{1}{\pi ir}\Big( \sinc(-\alpha r+s)-(-1)^{r}\delta_{\{s=0\} }\Big)
\end{eqnarray*}
while for $r=0$ and $s\neq 0$,
$$
c_\alpha (0,s)=\frac{-\alpha }{\pi is}\Big( \sinc(s)-\cos(\pi s)\Big)
=\frac{-\alpha }{\pi is}(-1)^{s}
$$
and $c_\alpha (0,0)=0$.

This way
\begin{eqnarray*}
{\mathcal D}_{\alpha }(a,b)(n)&=&\sum_{k_{1}\neq n}\sum_{k_{2}\in \bbbz }a_{k_{1}}b_{k_{2}}\frac{1}{\pi (n-k_{1})}\sinc((1-\alpha )n+k_{1}\alpha -k_{2})\\
&&-\sum_{k_{1}\neq n}a_{k_{1}}b_{n}\frac{1}{\pi (n-k_{1})}(-1)^{n-k_{1}}-\alpha \sum_{k_{2}\neq n}a_nb_{k_{2}}\frac{1}{\pi (n-k_{2})}(-1)^{n-k_{2}}\\
\end{eqnarray*}
Since $\sign(\eta +\alpha \xi )=\sign(\alpha )\sign (\xi +\frac{1}{\alpha }\eta )$, for $|\alpha |\geq 1$ we have
${\mathcal D}_{\alpha }(a,b)=\sign (\alpha ){\mathcal D}_{\frac{1}{\alpha }}(b,a)$ 
and so
$$
\sign (\alpha ){\mathcal D}_{\alpha }(a,b)(n)={\mathcal D}_{\frac{1}{\alpha }}(b,a)(n)
$$
$$
=\sum_{k_{1}\neq n}\sum_{k_{2}\in \bbbz }b_{k_{1}}a_{k_{2}}\frac{1}{\pi (n-k_{1})}
\sinc\big( (1-\frac{1}{\alpha })n+k_{1}\frac{1}{\alpha }-k_{2}\big)
$$
$$
-a_n\sum_{k_{1}\neq n}b_{k_{1}}\frac{1}{\pi (n-k_{1})}(-1)^{n-k_{1}}-\frac{1}{\alpha }b_n\sum_{k_{2}\neq n}a_{k_{2}}\frac{1}{\pi (n-k_{2})}(-1)^{n-k_{2}}
$$

Now, for $\alpha \in \mathbb Z\backslash \{ 0\}$ we define $\bar{a}$ by $\bar{a}_k=a_m$ if $k=\alpha m$ and zero otherwise and the same for $\bar{b}$. This way we have
$$
\sign (\alpha ){\mathcal D}_{\alpha }(\bar{a},\bar{b})(\alpha n)
$$
$$
=\sum_{k_{1}\neq n}\sum_{k_{2}\in \bbbz }b_{k_{1}}a_{k_{2}}\frac{1}{\pi \alpha (n-k_{1})}
\sinc\big( (1-\frac{1}{\alpha })\alpha n+\alpha k_{1}\frac{1}{\alpha }-\alpha k_{2}\big)
$$
$$
-a_{n}\sum_{k_{1}\neq n}b_{k_{1}}\frac{1}{\pi \alpha (n-k_{1})}(-1)^{\alpha (n-k_{1})}
-\frac{1}{\alpha }b_{n}\sum_{k_{2}\neq n}a_{k_{2}}\frac{1}{\pi \alpha (n-k_{2})}(-1)^{\alpha (n-k_{2})}
$$
$$
=\sum_{k_{1}\neq n}\sum_{k_{2}\in \bbbz }b_{k_{1}}a_{k_{2}}\frac{1}{\pi \alpha (n-k_{1})}
\sinc\big( (\alpha -1)n+k_{1}-\alpha k_{2}\big)
$$
$$
-\alpha^{-1}\tilde{a}_{n}{\mathcal H}(\tilde{b})(n)-\alpha^{-2}\tilde{b}_{n}{\mathcal H}(\tilde{a})(n)
$$
where $\tilde{a}_n=(-1)^{n}\bar{a}_n=(-1)^{\alpha m}a_{m}$ if $n=\alpha m$ and zero otherwise. The operator in the last line is 
the discrete linear Hilbert transform, 
defined by ${\mathcal H}(a)(n)=\frac{1}{\pi }\sum_{k\in \mathbb Z\backslash \{ 0\} }a_{n-k}\frac{1}{k}$.

Now since
$
(\alpha -1)n+k_{1}-\alpha k_{2}\in \mathbb Z
$, 
we have that 
$\sinc((\alpha -1)n+k_{1}-\alpha k_{2})$
is non zero only when $(\alpha -1)n+k_{1}$ is a multiple of $\alpha $ and
$\alpha k_{2}=(\alpha -1)n+k_{1}$. This leads to
$k_1
=n-\alpha (n-m)$ and
$k_{2}
=m$
and so
$$
\sign (\alpha ){\mathcal D}_{\alpha }(\bar{a},\bar{b})(\alpha n)
+\alpha^{-1}\tilde{a}_{n}H(\tilde{b})(n)+\alpha^{-2}\tilde{b}_{n}H(\tilde{a})(n)
$$
$$
=\sum_{m\neq n}b_{n-\alpha (n-m)}a_{m}\frac{1}{\pi \alpha^2 (n-m)}
=\frac{1}{\alpha^2\pi }\sum_{m\neq 0}b_{n-\alpha m}a_{n-m}\frac{1}{m}
$$
that is
$$
{\mathcal H}_{\alpha }(a,b)(n)
=\sign (\alpha )\alpha^{2}{\mathcal D}_{\alpha }(\bar{a},\bar{b})(\alpha n)
+\alpha \tilde{a}_{n}{\mathcal H}(\tilde{b})(n)+\tilde{b}_{n}{\mathcal H}(\tilde{a})(n)
$$

By using that the operators on the right are known to be bounded (the last two either by direct proof or by linear transference
from boundedness of the classical Hilbert tranform) and the inequality 
$
\sum_{n\in \mathbb Z}|{\mathcal D}_{\alpha }(\tilde{a},\tilde{b})(\alpha n)|^{p}
\leq \| {\mathcal D}_{\alpha }(\tilde{a},\tilde{b})\|_{l^{p}(\mathbb Z)}^{p}
$, we
finally deduce the stated boundedness of ${\mathcal H}_{\alpha }$ for $\alpha \in \mathbb Z\backslash \{ 0,1\}$.

\vskip 15pt
We end up with a remark about when the tranference results also extend to more general multilinear singular integrals.
The main use of the symbol properties is through the formula (\ref{tfoumulti})
$$
{\mathcal C}_{m}(f,g)\hate{}(\nu )
=\int_{\re }\hat{f}(\xi )\hat{g}(\nu -\xi )m(\nu ,\nu -\xi )d\xi
$$
which guarantees that if $f,g\in E_R$ then $C_m(f,g)\in E_{R'}$. For $x$-dependent multilinear operators like
$$
{\mathcal C}_{m}(f,g)(x)
=\int_{\re }\hat{f}(\xi )\hat{g}(\eta )m(x,\xi ,\eta )e^{2\pi ix(\xi+\eta )}d\xi d\eta
$$
this is also achieved if $m$ satisfies the property that the function $\widehat{m}^{x}(\nu, \xi, \eta ):=\int m(x,\xi ,\eta )e^{-2\pi ix(\nu-(\xi +\eta ))}dx $ (giving a distributional
sense to that integral if necessary) has compact support for all $\xi, \eta $ inside a compact set and such support is independent of $\xi ,\eta $. This way
$$
{\mathcal C}_{m}(f,g)\hate{}(\nu )
=\int_{\re }\hat{f}(\xi )\hat{g}(\eta )\int m(x,\xi ,\eta)e^{-2\pi ix(\nu-(\xi +\eta ))}dx d\xi
$$
and then the compact support ot the three functions in the integrand implies compact support for the Fourier transform of the operator itself. For such
operators a $T(1)$ type theorem like the one proved in \cite{T(1)} can be transferred to the discrete setting.

\bibliographystyle{plainnat}

\end{document}